\documentclass[12pt]{article}
\usepackage{mathrsfs}
\usepackage{latexsym,amsmath,amssymb,amsfonts,epsfig,graphicx,cite,psfrag}
\usepackage{eepic,color,colordvi,amscd,amsthm}
\usepackage{pstcol,pstricks,color}
\newtheorem{theorem}{Theorem}[section]

\newtheorem{definition}[theorem]{Definition}

 \numberwithin{equation}{section}
 \numberwithin{figure}{section}

\begin{document}

\begin{center}
{\Large\bf Pairs of Noncrossing Free Dyck Paths and

 Noncrossing
Partitions}
\end{center}

\vskip 2mm \centerline{William Y.C. Chen$^1$, Sabrina X.M. Pang$^2$,
Ellen X.Y. Qu$^3$, and Richard P. Stanley$^4$}

\begin{center}
$^{1, 2, 3}$Center for Combinatorics, LPMC-TJKLC\\
Nankai University, Tianjin 300071, P.R. China

$^4$Department of Mathematics\\
Massachusetts Institute of Technology, Cambridge, MA 02139, USA

\vskip 2mm
  $^1$chen@nankai.edu.cn,  $^2$pang@cfc.nankai.edu.cn,
 $^3$xiaoying@cfc.nankai.edu.cn, $^4$rstan@math.mit.edu
\end{center}



\begin{center}
{\bf Abstract}
\end{center}

Using the bijection between partitions and vacillating tableaux, we
establish a correspondence between pairs of noncrossing free Dyck
paths of length $2n$ and noncrossing partitions of $[2n+1]$ with
$n+1$ blocks. In terms of the number of up steps at odd positions,
we find a characterization of Dyck paths constructed from pairs of
noncrossing free Dyck paths by using the Labelle merging algorithm.

\vskip 3mm


\vskip 3mm


\section{Introduction}

We  use the bijection between vacillating tableaux and partitions to
establish a correspondence between pairs of noncrossing free Dyck
paths of length $2n$ and noncrossing partitions of $[2n+1]$ with
$n+1$ blocks. Recall that a {\it Dyck path} is a lattice path from
the origin to a point $(2n, 0)$ consisting of up steps $U=(1, 1)$
and down steps $D=(1, -1)$ that does not go below the $x$-axis.
Moreover, a lattice path from the origin to $(2n, 0)$ using the
steps $U$ and $D$ without the restriction on a Dyck path is called a
{\it
  free Dyck path}. Usually, a (free) Dyck path of length $2n$ is represented as a
sequence of $n$ $U$'s and $n$ $D$'s. A $k$-tuple $(P_1, P_2, \ldots
, P_k)$ of (free) Dyck paths from $(0,0)$ to $(2n,0)$ is called {\it
noncrossing} if each $P_i$ never goes below $P_{i+1}$ for $1\leq
i\leq k-1$.

A {\it partition} of a finite set $S$ is a collection $\pi=\{B_1,
B_2, \ldots, B_k\}$ of subsets of $S$ such that (i) $B_i \neq
\emptyset$ for each $i$; (ii) $B_i\cap B_j=\emptyset$ if $i\neq j$,
and (iii) $B_1 \cup B_2 \cup \cdots \cup B_k = S$. Each $B_i$ is
called a {\it block} of $\pi$. A {\it plane partition} is an array
$\delta=(\delta_{ij})_{i,j\geq1}$ of nonnegative integers such that
$\delta$ has finitely many nonzero entries and is weakly decreasing
in rows and columns. If $\sum\delta_{i,j}=n$, then we say that
$\delta$ is a plane partition of $n$ and write $|\delta|=n$. A {\it
part} of a plane partition $\delta=(\delta_{ij})$ is a positive
entry $\delta_{ij}>0$.  The shape of a plane partition $\delta$ is
the  integer partition $\lambda$ for which $\delta$  has $\lambda_i$
nonzero parts in the $i$-th row.

\begin{figure}
\begin{center}\setlength{\unitlength}{.3mm}
\begin{picture}(100,0)(0,0)
\multiput(0,0)(15,0){9}{\circle*{3}} \qbezier(0,0)(15,20)(30,0)
\qbezier(30,0)(45,20)(60,0)\qbezier(60,0)(82.5,25)(105,0)
\qbezier(45,0)(60,20)(75,0)\qbezier(15,0)(67.5,43)(120,0)
\put(-2,-14){\small$1$}\put(13,-14){\small$2$}
\put(28,-14){\small$3$} \put(43,-14){\small$4$}
\put(58,-14){\small$5$}\put(73,-14){\small$6$}
\put(88,-14){\small$7$}\put(103,-14){\small$8$}
\put(118,-14){\small$9$}
\end{picture}\vspace{-5pt}
\end{center}\caption{The standard
representation of $1358-29-46-7$.}\label{figure5}
\end{figure}
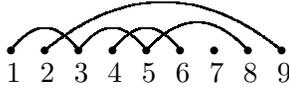

Let $[n]=\{1, 2,\ldots, n\}$. Given a partition $P$ of $[n]$, the
{\it standard representation} of the partition $P$ is a graph $G$ on
$[n]$ such that  a block $\{i_1, i_2,\ldots, i_k\}$ of $P$ written
in the increasing order $i_1< i_2< \cdots < i_k$ corresponds to a
path $(i_1, i_2, \ldots, i_k)$. For example, the standard
representation of $1358-29-46-7$ is illustrated in Figure
\ref{figure5}. Meanwhile, we may view the standard representation of
$P$ as a directed graph because each edge can always be considered
as an arc $(i, j)$ with $i < j$, and we say that $i$ is the {\it
left end point} and $j$ is the {\it right end point}. In this paper,
we use $\Pi_n$ to denote the set of partitions of $[n]$.

Let $k\geq 2$ and $P\in \Pi_n$. Define a {\it $k$-crossing} ({\it
$k$-nesting}) of $P$ as a set of $k$ arcs $(i_1, j_1), (i_2, j_2),
\ldots, (i_k, j_k)$ in the standard representation of $P$ such that
$i_1<i_2<\cdots<i_k<j_1<j_2<\cdots<j_k$
($i_1<i_2<\cdots<i_k<j_k<\cdots<j_2<j_1$). We use cr($P$) (ne($P$))
to denote the maximal $k$ such that $P$ has a $k$-crossing
($k$-nesting). In particular, a $2$-crossing ($2$-nesting)  is
called a {\it crossing (nesting)} for short.  The statistics cr and
ne were studied in \cite{CDDSY07} via a bijection between
vacillating tableaux and set-partitions. We will use this bijection
to connect  noncrossing free Dyck paths to noncrossing partitions.

 The paper is
organized as follows. In Section 2, we give a quick review of the
correspondence between $k$-tuples of noncrossing free Dyck paths and
plane partitions. As a consequence, we get the formula for the
number of pairs of noncrossing free Dyck paths of length $2n$. Then
we give an overview of the bijection between vacillating tableaux
and partitions as shown in \cite{CDDSY07}. We find a correspondence
between pairs of noncrossing free Dyck paths and vacillating
tableaux such that there is at most one row in each shape. These
vacillating tableaux allow us to construct the noncrossing
partitions. In Section 3,  we give a characterization of Dyck paths
obtained from pairs of noncrossing free Dyck paths by applying the
Labelle merging algorithm. 

\section{Pairs of Noncrossing Free Dyck Paths}
\label{sec2}

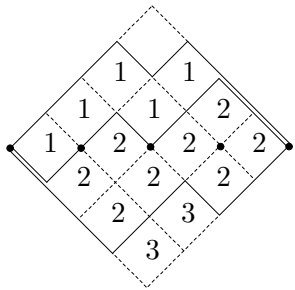
\begin{figure}
\begin{center}\setlength{\unitlength}{.07mm}
\begin{picture}(500,500)(0,-50)
\put(0,200){\line(1,1){202.6}} %
\put(-5,200){\line(1,-1){201}} %
\put(5,200){\line(1,-1){65}}  %
\put(70,135){\line(1,1){133}}  %
\multiput(70,265)(8,-8){33}{\line(1,-1){4}}%
\multiput(138,332)(8,-8){8}{\line(1,-1){4}} %
\multiput(273,332)(8,-8){24}{\line(1,-1){4}} %
\multiput(270,470)(8,-8){9}{\line(1,-1){4}} %
\multiput(206,271)(8,8){8}{\line(1,1){4}}  %
\multiput(134,69)(8,8){17}{\line(1,1){4}}  %
\multiput(335,140)(8,8){16}{\line(1,1){4}}  %
\multiput(262,-63)(8,8){17}{\line(1,1){4}}  %
\multiput(199,-4)(8,-8){8}{\line(1,-1){4}}  %
\put(195,0){\line(1,1){138}}  %
\put(332,138){\line(1,-1){66}}
\put(398,73){\line(1,1){135}}  %
\put(202,268){\line(1,-1){66}}
\put(267,202){\line(1,1){130}}  %
\put(397,332){\line(1,-1){130}} %
\put(203,403){\line(1,-1){68}} %
\multiput(208,408)(8,8){8}{\line(1,1){4}}%
\put(270,335){\line(1,1){67}}  %
\put(337,403){\line(1,-1){196}} %
\multiput(270,200)(8,-8){8}{\line(1,-1){4}} %
\put(0,200){\circle*{12}}  %
\put(135,200){\circle*{12}} %
\put(267,202){\circle*{12}} %
\put(400,203){\circle*{12}} %
\put(530,203){\circle*{12}} %
\put(64,192){\small$1$}   %
\put(128,258){\small$1$}%
\put(197,328){\small$1$} %
\put(195,192){\small$2$}
\put(260,258){\small$1$} %
\put(127,127){\small$2$}%
\put(192,60){\small$2$}
\put(259,127){\small$2$}%
\put(257,-10){\small$3$}%
\put(325,60){\small$3$}
\put(327,328){\small$1$} %
\put(327,192){\small$2$} %
\put(392,258){\small$2$} %
\put(392,127){\small$2$} %
\put(460,192){\small$2$} %
\end{picture}
\caption{A triple of noncrossing free Dyck paths and the
corresponding plane partition.}\label{figure1}
\end{center}
\end{figure}

We begin with the enumeration of $k$-tuples of noncrossing free Dyck
paths via the correspondence with plane partitions with bounded part
size. Given a $k$-tuple of noncrossing free Dyck paths
$(P_1,P_2,\ldots,P_k)$ with each $P_i$ of length $2n$, they must lie
in the region bounded by the paths $P_0$ consisting of $n$ up steps
followed by $n$ down steps and the path $P_{k+1}$ consisting of $n$
down steps followed by $n$ up steps. As illustrated in Figure
\ref{figure1}, we may obtain a plane partition by filling the areas
with $i$ in each square located in the region between the paths
$P_i$ and $P_{i+1}$ for $1\leq i\leq k$. Suppose that the resulting
plane partition $\pi$ is of shape $\lambda$. Then we see that
$\lambda_1\leq n$ and $\lambda'_1\leq n$ and the largest part of
$\pi$ does not exceed $k$. Since this correspondence is one-to-one,
the enumeration of $k$-tuples of noncrossing free Dyck paths can be
converted into the enumeration of plane partitions with bounded part
size.

Let $B(r,c,t)$ be the set of plane partitions with at most $r$ rows
and at most $c$ columns, and with the largest part at most $t$. It
is known   that
\begin{align}\label{bnnk}
\sum_{\pi\in B(n,n,k)}q^{|\pi|}=\frac
{[k+1][k+2]^2\cdots[k+n]^n[k+n+1]^{n-1}\cdots[k+2n-1]}
{[1][2]^2\cdots[n]^n[n+1]^{n-1}\cdots[2n-1]},
\end{align}
 where $[i]=1-q^i$; see, for example, \cite[Theorem 7.21.7]{Sta99}.
Let $F(n,k)$ denote the number of $k$-tuples of noncrossing free
Dyck paths  of length $2n$. Then one can deduce a formula for
$F(n,k)$ by setting $q=1$ in  \eqref{bnnk}. In particular,
\begin{equation}\label{e2.1}
 F(n,2)= \frac{(2n)!(2n+1)!}{(n!(n+1)!)^2}.
 \end{equation}
It has been shown by Callan that the above number also equals the
number of noncrossing partitions of $[2n+1]$ with $n+1$ blocks; see
Sloane \cite[Sequence A000891]{Slo}. Hence we are led to find a
bijection between the set of pairs of noncrossing free Dyck paths of
length $2n$ and the set of noncrossing partitions of $[2n+1]$ with
$n+1$ blocks.

Our bijection, denoted by  $\zeta$, consists of two steps. The first
step is to transform a pair of noncrossing free Dyck paths into a
vacillating tableau in which each shape has at most one row. Then we
use the bijection of Chen, Deng, Du, Stanley and Yan \cite{CDDSY07}
to construct the corresponding noncrossing partition. We now give a
brief review of the construction in \cite{CDDSY07}.  We assume that
the reader is familiar with the RSK algorithm, and we will use row
insertion as the basic operation.

\begin{definition}
A vacillating tableau $V_{\lambda}^{2n}$ of shape $\lambda$ and
length $2n$ is a sequence $(\lambda^0, \lambda^1, \ldots,
\lambda^{2n})$ of partitions such that (i) $\lambda^0=\emptyset$,
and $\lambda^{2n}=\lambda$, (ii) $\lambda^{2i+1}$ is obtained from
$\lambda^{2i}$ by doing nothing (i.e.,
$\lambda^{2i+1}=\lambda^{2i}$) or deleting a square, and (iii)
$\lambda^{2i}$ is obtained from $\lambda^{2i-1}$ by doing nothing or
adding a square.
\end{definition}

Given a partition $P$, let $E(P)$ denote the set of arcs in the
standard representation of $P$. To construct a vacillating tableau,
 we will derive a sequence $(T_0, T_1, T_2, \ldots,
T_{2n})$ of standard Young tableaux. Then the vacillating tableau is
just the sequence of the shapes of these tableaux. We work our way
backwards from $T_{2n}=\emptyset$ by determining the tableau
$T_{i-1}$ from $T_i$. We can construct the tableaux
$T_{2k-1},T_{2k-2}$ from $T_{2k}$ $(k\leq n)$ by the following
rules.

\begin{itemize}
\item[1.] Let $T_{2k-1}=T_{2k}$
if the integer $k$ does not appear in $T_{2k}$. Otherwise,
$T_{2k-1}$ is obtained from $T_{2k}$ by deleting the square occupied
by the element $k$.

\item[2.] $T_{2k-2}=T_{2k-1}$ if $E(P)$ does not have any arc of the
form $(i, k)$. Otherwise, there is a unique integer $i<k$ such that
$(i, k)\in E(P)$. Then   $T_{2k-2}$ is obtained from $T_{2k-1}$ by
row inserting the element $i$ into $T_{2k-1}$.
\end{itemize}
Let $\lambda^i$ be the shape of $T_i$ for $0\leq i\leq 2n$. Then
$(\lambda^0, \lambda^1,\lambda^2, \ldots, \lambda^{2n})$ is the
required vacillating tableau.

The inverse procedure can be described as follows. Given a
vacillating tableau $V=(\emptyset=\lambda^0, \lambda^1, \ldots,
\lambda^{2n}=\emptyset)$, we will recursively generate
 a sequence
$( T_0, T_1, \ldots,  T_{2n})$, where $T_i$ is a SYT (standard Young
tableau) of shape $\lambda^i$. Let $T_0$ be the empty SYT. Below are
the rules to construct $T_i$ from $ T_{i-1}$:

\begin{itemize}
\item[1.] If $\lambda^i=\lambda^{i-1}$, then  $T_i= T_{i-1}$.

\item[2.] If $\lambda^i\supset \lambda^{i-1}$, then $i=2k$ for some
integer $k\in[n]$. Determine the tableau $T_i$ such that $T_i$ is
obtained from $T_{i-1}$ by adding the integer $k$ in  the position
of $\lambda^i\setminus\lambda^{i-1}$.

\item[3.] If $\lambda^i\subset \lambda^{i-1}$,
then $i=2k-1$ for some integer $k\in[n]$. Set $T_i$ to be the unique
SYT (on a suitable alphabet) of shape $\lambda^i$ such that
$T_{i-1}$ is obtained from $T_i$ by row inserting some element $j$.
 Moreover, we record the arc $A_i=(j,k)$.
\end{itemize}
After the completion of the above procedure we are led to a set of
arcs  generated in Step 3. These arcs form a standard representation
of a partition  $P$ of $[n]$.

For example, given the vacillating tableau
\begin{center}
$(\emptyset, \emptyset, 1, 1, 11, 11, 111, 11, 11, 1, 2, 1, 1,
\emptyset, \emptyset)$,
\end{center}
the sequence $(A_i, T_i)$  is as follows:
\begin{center}
\begin{tabular}{c|lllllllllllllll}
$i$   & $0$         & $1$         & $2$ & $3$ & $4$  & $5$ & $6$ & $7$ & $8$ & $9$ & $10$ & $11$ & $12$ & $13$ & $14$ \\
\hline
$T_i$ & $\emptyset$ & $\emptyset$ & $1$ & $1$ & $1$ & $1$ & $1$ & $2$ & $2$ & $3$ & $35$ & $3$ & $3$ & $\emptyset$ & $\emptyset$ \\
      &   &   &   &   & $2$  &  $2$ &  $2$ &  $3$ & $3$ &  &  &  &  &  &  \\
      &   &   &   &   &   &   & $3$  &   &     &     &  &     &     &     &   \\
$A_i$ &   &   &   &   &   &   &   &  $(1, 4)$ &     &  $(2, 5)$   &     &$(5, 6)$  & & $(3, 7)$ &  \\
\end{tabular}
\end{center}
The corresponding  partition  is $P=14$-$256$-$37$.

The following result can be derived from \cite[Theorem 6]{CDDSY07}.

\begin{theorem}\label{theorem2.3} Let $P\in \Pi_n$ and
$(\emptyset=\lambda^0, \lambda^1, \ldots, \lambda^{2n}=\emptyset)$
be the corresponding vacillating tableau. Then cr($P$) is the
largest number of rows among $\lambda^i$, and ne($P$) is the largest
number of columns among $\lambda^i$.
\end{theorem}

We are now ready to describe the bijection $\zeta$. Let $(P,Q)$ be a
pair of noncrossing free Dyck paths, and let $P=p_1p_2\cdots p_{2n}$
and $Q=q_1q_2\cdots q_{2n}$. Based on $P$ and $Q$, we form the
sequence $(p_i, q_i)$, where $i=1, 2, \ldots, 2n$. The bijection
$\zeta$ consists of two phases. First,  we transform  $(P,Q)$ into a
vacillating tableau
\[
\mathcal{V}^{4n+2}_{\emptyset}=(\lambda^0,\lambda^1,\ldots,\lambda^{4n+2})\]
of empty shape, i.e. $\lambda^{4n+2}=\emptyset$, such that there is
at most one row in each $\lambda^i$ and there are a total number of
$n$ operations of adding a square in the process to obtain
$\lambda^{4n+2}$ from $\lambda^{0}$. Once a vacillating tableau is
constructed, we may turn the vacillating tableau into a partition by
the bijection in \cite{CDDSY07}.

For $1\leq i\leq 2n$, we have the following procedure to determine
$\lambda^{k}$ for $0\leq k\leq 4n+2$. Keep in mind that all the
involved tableaux have at most one row. Specifically, we have the
rules:
\begin{itemize}
\item[1.] $\lambda^0=\lambda^1=\lambda^{4n+2}=\emptyset$.
\item[2.] If $(p_i, q_i)=(U, U)$, then $\lambda^{2i}$ is
obtained from $\lambda^{2i-1}$ by adding one square, and
$\lambda^{2i+1}$ is obtained from $\lambda^{2i}$ by deleting one
square.

\item[3.] If $(p_i, q_i)=(U, D)$, then $\lambda^{2i}$ is obtained from $\lambda^{2i-1}$
by adding one square, and $\lambda^{2i+1} =\lambda^{2i}$.
\item[4.] If $(p_i, q_i)=(D, U)$, then
$\lambda^{2i}=\lambda^{2i-1}$, and $\lambda^{2i+1}$ is obtained from
$\lambda^{2i}$ by deleting one square.
\item[5.] If $(p_i, q_i)=(D, D)$, then  $\lambda^{2i-1}=\lambda^{2i}=
\lambda^{2i+1}$.
\end{itemize}

 Let $P_i=p_1p_2\cdots p_i$ and $Q_i=q_1q_2\cdots q_i$, and let
$|(P_i, Q_i)|_{(U, U)}$ denote the number of the pairs
$(p_j,q_j)=(U, U)$ in $(P_i, Q_i)$. Similarly, we can define $|(P_i,
Q_i)|_{(U, D)}$, $|(P_i, Q_i)|_{(D, U)}$ and $|(P_i, Q_i)|_{(D,
D)}$. Evidently, the number of $U$'s in $P_i$ equals $|(P_i,
Q_i)|_{(U, U)}+|(P_i, Q_i)|_{(U, D)}$, and the number of $U$'s in
$Q_i$ equals $|(P_i, Q_i)|_{(U, U)}
 +|(P_i, Q_i)|_{(D, U)}$.
Since $P$ and $Q$ are noncrossing, the number of $U$'s in $P_i$ is
 not less than that in $Q_i$. It follows that
\[
 |(P_i, Q_i)|_{(U, U)}+|(P_i, Q_i)|_{(U, D)}\geq |(P_i,
Q_i)|_{(U, U)}
 +|(P_i, Q_i)|_{(D, U)}.
\]
Hence we find that
 \begin{equation} \label{pq}
 |(P_i, Q_i)|_{(U, D)}\geq |(P_i, Q_i)|_{(D, U)}.
\end{equation}

To show that $(\lambda^0, \lambda^1, \lambda^2, \ldots,
\lambda^{4n+2})$ is a valid vacillating tableau, we should justify
that the above constructions are feasible. Clearly, the items 1, 2,
3 and 5 are well defined. So we may restrict our attention  to item
4, in which case $(p_i, q_i)=(D, U)$. We aim to show that
$\lambda^{2i}=\lambda^{2i-1}$ is not an empty shape. Combining
(\ref{pq}) and the relations $ |(P_i, Q_i)|_{(D, U)}=|(P_{i-1},
Q_{i-1})|_{(D, U)}+1$ and $ |(P_i, Q_i)|_{(U, D)}=|(P_{i-1},
Q_{i-1})|_{(U, D)},$ we see that
\[
|(P_{i-1}, Q_{i-1})|_{(U, D)}> |(P_{i-1}, Q_{i-1})|_{(D, U)}.\] From
the construction of $(\lambda^0, \lambda^1, \ldots, \lambda^{2i-1})$
it can be seen that $\lambda^{2i-1}\not=\emptyset$. Thus we have
reached the conclusion that $(\lambda^0, \lambda^1, \ldots,
\lambda^{4n+2})$ is a vacillating tableau.

It is easily seen from Theorem \ref{theorem2.3} that the partition
corresponding to the above  vacillating tableau, denoted by $R$,
 is
noncrossing because $\lambda^i$ contains at most one row for any
$i$. It remains to show that the resulting partition $R$ contains
exactly $n+1$ blocks. Since the free Dyck path $P$ is of length
$2n$, there are $n$ left end points in the standard representation
of $R$. This implies that there are $n$ arcs in the standard
representation of $R$. On the other hand, $R$ contains $2n+1$
elements, since the vacillating tableau is of length $4n+2$. So we
may deduce that there are $n+1$ blocks in $R$. It is not difficult
to see that the above procedure is reversible. Therefore, we have
established the following result.

\begin{theorem}\label{theorem2.4}
The above map $\zeta$ is a bijection between the set of pairs of
noncrossing free Dyck paths of length $2n$ and the set of
noncrossing partitions of $[2n+1]$ with $n+1$ blocks.
\end{theorem}

Figure \ref{figure2} is an illustration of the bijection $\zeta$.

\begin{figure}
\begin{center}\setlength{\unitlength}{.25mm}
\begin{picture}(150,10)(0,0)
\put(0,0){\circle*{2}}
\put(0,0){\line(1,1){20}}\put(20,20){\circle*{2}}
\put(20,20){\line(1,1){20}}\put(40,40){\circle*{2}}
\put(40,40){\line(1,-1){20}}\put(60,20){\circle*{2}}
\put(60,20){\line(1,-1){20}}\put(80,0){\circle*{2}}
\put(80,0){\line(1,-1){20}}\put(100,-20){\circle*{2}}
\put(100,-20){\line(1,1){20}}\put(120,0){\circle*{2}}
\put(120,0){\line(1,1){20}}\put(140,20){\circle*{2}}
\put(140,20){\line(1,-1){20}}\put(160,0){\circle*{2}}
\put(0,-5){\circle*{2}}\put(0,-5){\line(1,-1){20}}
\put(20,-25){\circle*{2}}\put(20,-25){\line(1,1){20}}
\put(40,-5){\circle*{2}}\put(40,-5){\line(1,1){20}}
\put(60,15){\circle*{2}}\put(60,15){\line(1,-1){20}}
\put(80,-5){\circle*{2}}\put(80,-5){\line(1,-1){20}}
\put(100,-25){\circle*{2}}\put(100,-25){\line(1,1){20}}
\put(120,-5){\circle*{2}}\put(120,-5){\line(1,-1){20}}
\put(140,-25){\circle*{2}}\put(140,-25){\line(1,1){20}}
\put(160,-5){\circle*{2}}
\put(70,-40){$\Updownarrow$}
\put(-150,-65){\small$T_0$}\put(-125,-65){\small$T_1$}
\put(-100,-65){\small$T_2$} \put(-75,-65){\small$T_3$}
\put(-50,-65){\small$T_4$} \put(-25,-65){\small$T_5$}
\put(0,-65){\small$T_6$} \put(25,-65){\small$T_7$}
\put(50,-65){\small$T_8$}\put(75,-65){\small$T_9$}
\put(100,-65){\small$T_{10}$} \put(125,-65){\small$T_{11}$}
\put(150,-65){\small$T_{12}$} \put(175,-65){\small$T_{13}$}
\put(200,-65){\small$T_{14}$} \put(225,-65){\small$T_{15}$}
\put(250,-65){\small$T_{16}$} \put(275,-65){\small$T_{17}$}
\put(300,-65){\small$T_{18}$}
\put(-150,-85){$\emptyset$}\put(-125,-85){$\emptyset$}
\put(-100,-87){$\Box$}\put(-75,-87){$\Box$}
\put(-54,-87){$\Box$}\put(-45,-87){$\Box$} \put(-25,-87){$\Box$}
\put(0,-87){$\Box$}\put(25,-85){$\emptyset$}
\put(50,-85){$\emptyset$} \put(75,-85){$\emptyset$}
\put(100,-85){$\emptyset$} \put(125,-85){$\emptyset$}
\put(150,-87){$\Box$}\put(175,-85){$\emptyset$}
\put(200,-87){$\Box$} \put(225,-87){$\Box$} \put(250,-87){$\Box$}
\put(275,-85){$\emptyset$}\put(300,-85){$\emptyset$}
\put(70,-110){$\Updownarrow$}
\put(-10,-140){\circle*{2}}\put(10,-140){\circle*{2}}
\put(30,-140){\circle*{2}} \put(50,-140){\circle*{2}}
\put(70,-140){\circle*{2}} \put(90,-140){\circle*{2}}
\put(110,-140){\circle*{2}} \put(130,-140){\circle*{2}}
\put(150,-140){\circle*{2}} \put(-12,-154){\small$1$}
\put(8,-154){\small$2$} \put(28,-154){\small$3$}
\put(48,-154){\small$4$} \put(68,-154){\small$5$}
\put(88,-154){\small$6$}\put(108,-154){\small$7$}
\put(128,-154){\small$8$}\put(148,-154){\small$9$}
\qbezier(-10,-140)(20,-100)(50,-140)
\qbezier(10,-140)(20,-120)(30,-140)
\qbezier(90,-140)(100,-120)(110,-140)
\qbezier(110,-140)(130,-110)(150,-140)
\end{picture}\vspace{105pt}
\caption{A pair of noncrossing free Dyck paths and its corresponding
partition.} \label{figure2}
\end{center}
\end{figure}

We remark that the bijection $\zeta$ can be described in a simpler
 manner. For a given pair $(P,Q)$ of noncrossing free Dyck paths,
let $P=p_1p_2\cdots p_{2n}$ and $Q=q_1q_2\cdots q_{2n}$, and let
$l_i$ (resp. $r_{i}$) denote the left-degree (resp. the
right-degree) of vertex $i$ in the standard representation of the
partition corresponding to $(P,Q)$, i.e., the number of vertices $j$
with $j<i$ (resp. $j>i$) connected to $i$. First, set
$l_1=r_{2n+1}=0$. Then the pair $(r_i,l_{i+1})$ is determined by the
following rules for $1\leq i\leq 2n$.
\begin{itemize}
\item[1.] If $(p_i, q_i)=(U, U)$, then set $(r_i, l_{i+1})=(1, 1)$.
\item[2.] If $(p_i, q_i)=(U, D)$, then set $(r_i, l_{i+1})=(1, 0)$.
\item[3.] If $(p_i, q_i)=(D, U)$, then set $(r_i, l_{i+1})=(0, 1)$.
\item[4.] If $(p_i, q_i)=(D, D)$, then set $(r_i, l_{i+1})=(0, 0)$.
\end{itemize}

Observe that the degree sequences $(l_1, l_2, \ldots, l_{2n+1})$ and
$(r_1, r_2, \ldots, r_{2n+1})$ consisting of only zeros and ones. As
shown in Figure \ref{figure3}, we may use half arcs (intuitively
called left half arcs and right half arcs)  to represent the left
and right degrees.  We have the following unique way to pair up half
arcs in order to form a noncrossing partition. At each step we
always try to find the leftmost left half arc and pair it up with
the nearest right half arc on its left. Iterate this procedure until
all the half arcs are paired. Finally, we obtain the standard
representation of the desired partition, say $R$.

\begin{figure}
\begin{center}\setlength{\unitlength}{.25mm}
\begin{picture}(150,10)(0,0)
\put(-10,0){\circle*{2}}
\put(-10,0){\line(1,1){20}}\put(10,20){\circle*{2}}
\put(10,20){\line(1,1){20}}\put(30,40){\circle*{2}}
\put(30,40){\line(1,-1){20}}\put(50,20){\circle*{2}}
\put(50,20){\line(1,-1){20}}\put(70,0){\circle*{2}}
\put(70,0){\line(1,-1){20}}\put(90,-20){\circle*{2}}
\put(90,-20){\line(1,1){20}}\put(110,0){\circle*{2}}
\put(110,0){\line(1,1){20}}\put(130,20){\circle*{2}}
\put(130,20){\line(1,-1){20}}\put(150,0){\circle*{2}}
\put(-10,-5){\circle*{2}}\put(-10,-5){\line(1,-1){20}}
\put(10,-25){\circle*{2}}\put(10,-25){\line(1,1){20}}
\put(30,-5){\circle*{2}}\put(30,-5){\line(1,1){20}}
\put(50,15){\circle*{2}}\put(50,15){\line(1,-1){20}}
\put(70,-5){\circle*{2}}\put(70,-5){\line(1,-1){20}}
\put(90,-25){\circle*{2}}\put(90,-25){\line(1,1){20}}
\put(110,-5){\circle*{2}}\put(110,-5){\line(1,-1){20}}
\put(130,-25){\circle*{2}}\put(130,-25){\line(1,1){20}}
\put(150,-5){\circle*{2}}
\put(70,-40){$\Updownarrow$}
\put(-10,-70){\circle*{2}}\put(10,-70){\circle*{2}}
\put(30,-70){\circle*{2}} \put(50,-70){\circle*{2}}
\put(70,-70){\circle*{2}} \put(90,-70){\circle*{2}}
\put(110,-70){\circle*{2}} \put(130,-70){\circle*{2}}
\put(150,-70){\circle*{2}} \put(-12,-84){\small$1$}
\put(8,-84){\small$2$} \put(28,-84){\small$3$}
\put(48,-84){\small$4$} \put(68,-84){\small$5$}
\put(88,-84){\small$6$}\put(108,-84){\small$7$}
\put(128,-84){\small$8$}\put(148,-84){\small$9$}
\qbezier(-10,-70)(-9,-64)(-5,-60) \qbezier(10,-70)(11,-64)(15,-60)
\qbezier(30,-70)(29,-64)(25,-60)\qbezier(50,-70)(49,-64)(45,-60)
\qbezier(90,-70)(91,-64)(95,-60)\qbezier(110,-70)(109,-64)(105,-60)
\qbezier(110,-70)(111,-64)(115,-60)\qbezier(150,-70)(149,-64)(145,-60)
\put(70,-110){$\Updownarrow$}
\put(-10,-125){\circle*{2}}\put(10,-125){\circle*{2}}
\put(30,-125){\circle*{2}} \put(50,-125){\circle*{2}}
\put(70,-125){\circle*{2}} \put(90,-125){\circle*{2}}
\put(110,-125){\circle*{2}} \put(130,-125){\circle*{2}}
\put(150,-125){\circle*{2}} \put(-12,-137){\small$1$}
\put(8,-139){\small$2$} \put(28,-139){\small$3$}
\put(48,-139){\small$4$} \put(68,-139){\small$5$}
\put(88,-139){\small$6$}\put(108,-139){\small$7$}
\put(128,-139){\small$8$}\put(148,-139){\small$9$}
\qbezier(-10,-125)(20,-85)(50,-125)
\qbezier(10,-125)(20,-105)(30,-125)
\qbezier(90,-125)(100,-105)(110,-125)
\qbezier(110,-125)(130,-95)(150,-125)
\end{picture}
\end{center}\vspace{95pt}
\caption{A pair of noncrossing free Dyck paths and the corresponding
partition.}\label{figure3}
\end{figure}

It is necessary to show that the  resulting partition  $R$ contains
$2n+1$ elements and  $n+1$ blocks and is noncrossing. First, from
the sequence $(p_1, q_1), (p_2, q_2), \ldots, (p_{2n}, q_{2n})$ we
determine that degrees $(r_1, l_2), (r_2, l_3), \ldots, (r_{2n},
l_{2n+1})$. Clearly,  the underlying set of the partition $R$ is
$[2n+1]$. From the definition of $(r_i, l_{i+1})$, one sees that the
total number of half arcs equals the number of $U$'s in $P$ and $Q$.
Since $P$ and $Q$ are free Dyck paths of length $2n$, the total
number of $U$'s in $P$ and $Q$ equals $2n$. It follows that there
are $n$ arcs in the standard representation of $R$. Moreover, $R$ is
noncrossing since at each step the leftmost left half arc is paired
up with the nearest right half arc on its left, and this operation
does not cause any crossing. It is not hard to check that the above
procedure is reversible.  Figure \ref{figure3} gives an illustration
of this procedure.

\section{The Labelle  Merging Algorithm}

In this section, we establish a correspondence between pairs of
noncrossing free Dyck paths of length $2n$ and Dyck paths of length
$4n+2$ with $n+1$ up steps at odd positions. According to a theorem
of Sulanke \cite{Sul05},
 the number of Dyck paths  of length $2n$ with $k$
 up steps at odd positions equals the Narayana number
 $N(n,k) = \frac 1n
                 \binom{n}{k-1} \binom nk$. This implies the formula
                 \eqref{e2.1}.

 Labelle \cite{Lab93} gives an algorithm to merge a pair
of noncrossing free Dyck paths into a $2$-Motzkin path. It is
realized that one can further transform a $2$-Motzkin path into a
Dyck path by using a bijection due to Delest and Viennot
\cite{Del-Vie84}. We will present an equivalent algorithm that
directly transforms a pair of noncrossing free Dyck paths into a
single Dyck path, and will call it the Labelle merging algorithm.

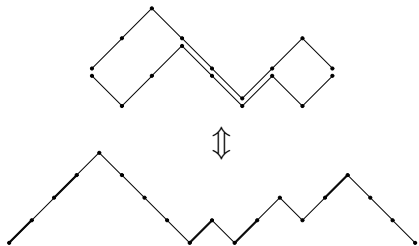
\begin{figure}
\begin{center}\vspace{20pt}
\setlength{\unitlength}{.2mm}
\begin{picture}(150,20)(0,0)
\put(-10,0){\circle*{2}}
\put(-10,0){\line(1,1){20}}\put(10,20){\circle*{2}}
\put(10,20){\line(1,1){20}}\put(30,40){\circle*{2}}
\put(30,40){\line(1,-1){20}}\put(50,20){\circle*{2}}
\put(50,20){\line(1,-1){20}}\put(70,0){\circle*{2}}
\put(70,0){\line(1,-1){20}}\put(90,-20){\circle*{2}}
\put(90,-20){\line(1,1){20}}\put(110,0){\circle*{2}}
\put(110,0){\line(1,1){20}}\put(130,20){\circle*{2}}
\put(130,20){\line(1,-1){20}}\put(150,0){\circle*{2}}
\put(-10,-5){\circle*{2}}\put(-10,-5){\line(1,-1){20}}
\put(10,-25){\circle*{2}}\put(10,-25){\line(1,1){20}}
\put(30,-5){\circle*{2}}\put(30,-5){\line(1,1){20}}
\put(50,15){\circle*{2}}\put(50,15){\line(1,-1){20}}
\put(70,-5){\circle*{2}}\put(70,-5){\line(1,-1){20}}
\put(90,-25){\circle*{2}}\put(90,-25){\line(1,1){20}}
\put(110,-5){\circle*{2}}\put(110,-5){\line(1,-1){20}}
\put(130,-25){\circle*{2}}\put(130,-25){\line(1,1){20}}
\put(150,-5){\circle*{2}} \put(70,-55){$\Updownarrow$}
\end{picture}
\end{center}
\begin{center}\vspace{30pt}
\setlength{\unitlength}{.2mm}
\begin{picture}(230,20)(0,0)
\linethickness{.5pt} \put(-25,-15){\circle*{2}}
\put(-25.5,-15){\line(1,1){15}}\put(-24.5,-15){\line(1,1){15}}\put(-25,-15){\line(1,1){15}}
 \put(-10,0){\circle*{2}}
\put(-10,0){\line(1,1){15}} \put(5,15){\circle*{2}}
\put(4.5,15){\line(1,1){15}}\put(5.5,15){\line(1,1){15}}\put(5,15){\line(1,1){15}}
 \put(20,30){\circle*{2}}
\put(20,30){\line(1,1){15}} \put(35,45){\circle*{2}}
\put(35,45){\line(1,-1){15}}\put(50,30){\circle*{2}}
\put(50,30){\line(1,-1){15}}\put(65,15){\circle*{2}}
\put(65,15){\line(1,-1){15}} \put(80,0){\circle*{2}}
\put(80,0){\line(1,-1){15}}\put(95,-15){\circle*{2}}
\put(94.5,-15){\line(1,1){15}}\put(95.5,-15){\line(1,1){15}}\put(95,-15){\line(1,1){15}}
\put(110,0){\circle*{2}}
\put(110,0){\line(1,-1){15}}\put(125,-15){\circle*{2}}
\put(124.5,-15){\line(1,1){15}}\put(125.5,-15){\line(1,1){15}}\put(125,-15){\line(1,1){15}}
\put(140,0){\circle*{2}}
\put(140,0){\line(1,1){15}}\put(155,15){\circle*{2}}
\put(155,15){\line(1,-1){15}} \put(170,0){\circle*{2}}
\put(170,0){\line(1,1){15}} \put(185,15){\circle*{2}}
\put(184.5,15){\line(1,1){15}}\put(185.5,15){\line(1,1){15}}\put(185,15){\line(1,1){15}}
\put(200,30){\circle*{2}}
\put(200,30){\line(1,-1){15}}\put(215,15){\circle*{2}}
\put(215,15){\line(1,-1){15}} \put(230,0){\circle*{2}}
\put(230,0){\line(1,-1){15}}\put(245,-15){\circle*{2}}
\end{picture}
\end{center}\vspace{5pt}
\caption{The Labelle merging algorithm.}\label{figure4}
\end{figure}

Let $P=p_1p_2\cdots p_{2n}$ and $Q=q_1q_2\cdots q_{2n}$ be a pair of
noncrossing free Dyck paths of length $2n$. Let $Q'=
    q_1'q_2' \cdots q_{2n}'$, where we
    define $U'=D$ and $D'=U$.
    Then we merge $P$ and $Q$ into a Dyck path
    \[ Up_1q_1'p_2q_2'
\cdots p_{2n}q_{2n}'D.\]
 The following theorem gives a
characterization of the Dyck paths corresponding to pairs of
noncrossing free Dyck paths.

\begin{theorem}\label{theorem3.1}
The Labelle merging algorithm is a bijection between noncrossing
free Dyck paths of length $2n$ and Dyck paths of length $4n+2$ with
$n+1$ up steps at odd positions.
\end{theorem}

The verification of the above statement is omitted.
 Figure \ref{figure4} is an illustration of the Labelle  merging
algorithm, where the thick lines represent up steps at odd
positions.

\vskip 3mm \noindent {\bf Acknowledgments.} The authors would like
to thank Martin Rubey and the referee for valuable comments leading
to an improvement of an earlier version. This work was supported by
the 973 Project, the PCSIRT Project of the Ministry of Education,
the Ministry of Science and Technology, the National Science
Foundation of China, and the US NSF grants DMS-9988459 and
DMS-0604423.

\end{document}